\documentclass[12]{amsart}
\usepackage{amsmath,amssymb,amsthm,color,enumerate,comment,centernot,enumitem,url,cite}
\usepackage{graphicx,relsize,bm}
\usepackage{mathtools}
\usepackage{array}

\makeatletter
\newcommand{\tpmod}[1]{{\@displayfalse\pmod{#1}}}
\makeatother

\newtheorem{thm}{Theorem}[section]

\theoremstyle{remark}

\theoremstyle{definition}

\newtheorem{rem}[thm]{Remark}

\theoremstyle{THM}

\newcommand{\abs}[1]{\left|{#1}\right|}

\def\HH {{\mathcal H}}
\def\FF {{\mathcal F}}
\def\GG {{\mathcal G}}

\def\R {{\mathbb R}}
\def\Z {{\mathbb Z}}

\def\Q {{\mathbb Q}}

\def\GG {{\mathcal G}}

\def\F {{\mathbb F}}

\def\Z {{\mathbb Z}}
\def\Q {{\mathbb Q}}

\def\Gal{{{\rm Gal} }}

\makeatletter
\@namedef{subjclassname@2020}{%
  \textup{2020} Mathematics Subject Classification}
\makeatother

\def\red#1 {\textcolor{red}{#1 }}
\def\blue#1 {\textcolor{blue}{#1 }}

\numberwithin{equation}{section}

\def\Z {{\mathbb Z}}

\newcommand{\Mod}[1]{\ (\mathrm{mod}\enspace #1)}

\begin{document}

\title[Monogenic Quartics and Their Galois Groups]{Monogenic Quartic Polynomials\\ and Their Galois Groups}

\author{Joshua Harrington}
\address{Department of Mathematics, Cedar Crest College, Allentown, Pennsylvania, USA}
\email[Joshua Harrington]{Joshua.Harrington@cedarcrest.edu}

\author{Lenny Jones}
\address{Professor Emeritus, Department of Mathematics, Shippensburg University, Shippensburg, Pennsylvania 17257, USA}
\email[Lenny~Jones]{doctorlennyjones@gmail.com}

\date{\today}

\begin{abstract}
A monic polynomial $f(x)\in \Z[x]$ of degree $N$ is called \emph{monogenic} if $f(x)$ is irreducible over $\Q$ and $\{1,\theta,\theta^2,\ldots ,\theta^{N-1}\}$ is a basis for the ring of integers of $\Q(\theta)$, where $f(\theta)=0$. In this article, we use the classification of the Galois groups of quartic polynomials, due to Kappe and Warren, to investigate the existence of infinite collections of monogenic quartic polynomials having a prescribed Galois group, such that each member of the collection generates a distinct quartic field. With the exception of the cyclic case, we provide such an infinite single-parameter collection for each possible Galois group. We believe these examples are new, and we provide evidence to support this belief by showing that they are distinct from other infinite collections in the current literature. Finally, we devote a separate section to a discussion concerning, what we believe to be, the still-unresolved cyclic case.
\end{abstract}


\maketitle
\section{Introduction}\label{Section:Intro}

  We say that a monic polynomial $f(x)\in \Z[x]$ is \emph{monogenic} if $f(x)$ is irreducible over $\Q$ and $\{1,\theta,\theta^2,\ldots ,\theta^{\deg{f}-1}\}$ is an integral basis for the ring of integers $\Z_K$ of $K=\Q(\theta)$, where $f(\theta)=0$; that is, $\Z_K=\Z[\theta]$. We let $\Delta(f)$ and $\Delta(K)$ denote the discriminants over $\Q$, respectively, of $f(x)$ and the number field $K$. When $f(x)$ is irreducible over $\Q$, we have \cite{Cohen}
\begin{equation} \label{Eq:Dis-Dis}
\Delta(f)=\left[\Z_K:\Z[\theta]\right]^2\Delta(K),
\end{equation}
so that, in this situation, $f(x)$ is monogenic if and only if
  $\Delta(f)=\Delta(K)$. We also say that any number field $K$ is \emph{monogenic} if there exists a power basis for $\Z_K$.
We point out to the reader that, while the monogenicity of $f(x)$ implies the monogenicity of $K=\Q(\theta)$, where $f(\theta)=0$,
 the converse is not necessarily true. For example, let $f(x)=x^2-5$ and $K=\Q(\theta)$, where $\theta=\sqrt{5}$. Then, easy calculations show that
$\Delta(f)=20$ and $\Delta(K)=5$. Thus, $f(x)$ is not monogenic, but nevertheless, $K$ is monogenic since $\{1,(\theta+1)/2\}$ is a power basis for $\Z_K$. Note then that $g(x)=x^2-x-1$, the minimal polynomial for $(\theta+1)/2$, is monogenic.

It is the goal of this article to determine new families of monogenic quartic polynomials that have a prescribed Galois group over $\Q$. A ``family" here means an infinite collection of polynomials, with each polynomial in the family generating a distinct field. For quartic polynomials $f(x)$ that are irreducible over $\Q$, there are five possibilities for the Galois group $\Gal(f)$ over $\Q$, which are given in Table \ref{T:1} using the standard ``4TX" notation \cite{BM,Dokchitser}, along with more familiar names.
 \begin{table}[h]
 \begin{center}
\begin{tabular}{c|ccccccc}
 X &  1 & 2 & 3 & 4 & 5  \\ [2pt]  \hline \\[-8pt]
 4TX & $C_4$ & $C_2\times C_2$ & $D_4$ & $A_4$ & $S_4$
 \end{tabular}
\end{center}
\caption{The transitive groups of degree 4}
 \label{T:1}
\end{table}
With the exception of the case of 4T1, we provide a single-parameter family of monogenic quartic polynomials for each group in Table \ref{T:1}, while also providing evidence in a separate section (see Section \ref{Section:Known}) to support our belief that these families are actually new. Finally, since we have not been able to produce a family of cyclic quartic polynomials, we devote a separate section to a discussion concerning, what we believe to be, the still-unresolved  case of 4T1 (see Section \ref{Section:4T1}). The new results are presented in the following main theorem.
\begin{thm}\label{Thm:Main1}\text{}
  \begin{itemize}
  \item \label{I:2} {\bf X=2}\\
  Let $t\in \Z$ and let $f_t(x):=x^4+4tx^2+1$. Then
  \begin{enumerate}
 \item \label{X=2:I1} $f_t(x)$ is irreducible and $\Gal(f_t)\simeq$ {\rm 4T2},
 \item \label{X=2:I2} $f_t(x)$ is monogenic if and only if $(2t-1)(2t+1)$ is squarefree,
 \item \label{X=2:I3}
 $\FF_2:=\{f_t(x): 4t^2-1 \mbox{ is squarefree}\}$\\
 is an infinite family of monogenic {\rm 4T2}-quartics.
  \end{enumerate}
  \item \label{I:3} {\bf X=3}\\
  Let $t\in \Z$ and let $f_t(x):=x^4+24tx^3+(12t+4)x^2+4x+1$. Then
  \begin{enumerate}
 \item \label{X=3:I1} $f_t(x)$ is irreducible and $\Gal(f_t)\simeq$ {\rm 4T3},
 \item \label{X=3:I2} $f_t(x)$ is monogenic if and only if $(6t-1)(6t+1)$ is squarefree,
 \item \label{X=3:I3}
 $\FF_3:=\{f_t(x): 36t^2-1 \mbox{ is squarefree}\}$\\
 is an infinite family of monogenic {\rm 4T3}-quartics.
  \end{enumerate}
  \item \label{I:4} {\bf X=4}\\
  Let $t\in \Z$ and let $f_t(x):=x^4+2x^3+2x^2+4tx+36t^2-16t+2$. Then
  \begin{enumerate}
 \item \label{X=4:I1} $f_t(x)$ is irreducible and $\Gal(f_t)\simeq$ {\rm 4T4},
 \item \label{X=4:I2} $f_t(x)$ is monogenic if and only if $(4t-1)(108t^2-54t+7)$ is squarefree,
 \item \label{X=4:I3}
 $\FF_4:=\{f_t(x): (4t-1)(108t^2-54t+7) \mbox{ is squarefree}\}$\\
 is an infinite family of monogenic {\rm 4T4}-quartics.
  \end{enumerate}
 \item \label{I:5} {\bf X=5}\\
  Let $t\in \Z$, and let $f_t(x):=x^4-2x^3-2x^2+6x+4t-2$. Then
  \begin{enumerate}
 \item \label{X=5:I1} $f_t(x)$ is irreducible and $\Gal(f_t)\simeq$ {\rm 4T5},
 \item \label{X=5:I2} $f_t(x)$ is monogenic if and only if $4t+1$, $4t-7$ and $64t+13$ are squarefree,
 \item \label{X=5:I3}
 $\FF_5:=\{f_t(x):  4t+1,\ 4t-7 \mbox{ and } 64t+13 \mbox{ are squarefree}\}$\\
 is an infinite family of monogenic {\rm 4T5}-quartics.
  \end{enumerate}
\end{itemize}
\end{thm}
To the best of our knowledge, the families given in Theorem \ref{Thm:Main1} are new, and in Section \ref{Section:Known}, we show that these families contain no overlap with all explicit families that we found in the existing literature.

\section{Preliminaries}\label{Section:Prelim}
The following theorem is due to Kappe and Warren \cite{KW}.
\begin{thm}{\rm \cite{KW}}\label{Thm:KW1}
Let $f(x)=x^4+ax^3+bx^2+cx+d\in \Z[x]$ be irreducible over $\Q$. Let \[r(x)=x^3-bx^2+(ac-4d)x-(a^2d-4bd+c^2)\] with splitting field $L$. Then $\Gal(f)\simeq$
\begin{enumerate}
 \item \label{GI:1} {\rm 4T1} if and only if $r(x)$ has exactly one root $s\in \Q$ and
 \begin{equation}\label{Eq:gr}
 g(x):=(x^2-sx+d)(x^2+ax+(b-s))
 \end{equation} splits over $L$;
 \item \label{GI:2} {\rm 4T2} if and only if $r(x)$ splits into linear factors over $\Q$;
 \item \label{GI:3} {\rm 4T3} if and only if $r(x)$ has exactly one root $s\in \Q$ and $g(x)$, as defined in \eqref{Eq:gr}, does not split over $L$;
 \item \label{GI:4} {\rm 4T4} if and only if $r(x)$ is irreducible over $\Q$ and $\Delta(f)$ is a square in $\Z$;
  \item \label{GI:5} {\rm 4T5} if and only if $r(x)$ is irreducible over $\Q$ and $\Delta(f)$ is not a square in $\Z$;
 \end{enumerate}
 \end{thm}
 \begin{rem}
   The polynomial $r(x)$ in Theorem \ref{Thm:KW1} is known as \emph{the cubic resolvent of $f(x)$}.
 \end{rem}

The following theorem, known as \emph{Dedekind's Index Criterion}, or simply \emph{Dedekind's Criterion} if the context is clear, is a standard tool used in determining the monogenicity of an irreducible polynomial.
\begin{thm}[Dedekind \cite{Cohen}]\label{Thm:Dedekind}
Let $K=\Q(\theta)$ be a number field, $T(x)\in \Z[x]$ the monic minimal polynomial of $\theta$, and $\Z_K$ the ring of integers of $K$. Let $q$ be a prime number and let $\overline{ * }$ denote reduction of $*$ modulo $q$ (in $\Z$, $\Z[x]$ or $\Z[\theta]$). Let
\[\overline{T}(x)=\prod_{i=1}^k\overline{\tau_i}(x)^{e_i}\]
be the factorization of $T(x)$ modulo $q$ in $\F_q[x]$, and set
\[h_1(x)=\prod_{i=1}^k\tau_i(x),\]
where the $\tau_i(x)\in \Z[x]$ are arbitrary monic lifts of the $\overline{\tau_i}(x)$. Let $h_2(x)\in \Z[x]$ be a monic lift of $\overline{T}(x)/\overline{h_1}(x)$ and set
\[F(x)=\dfrac{h_1(x)h_2(x)-T(x)}{q}\in \Z[x].\]
Then
\[\left[\Z_K:\Z[\theta]\right]\not \equiv 0 \pmod{q} \Longleftrightarrow \gcd\left(\overline{h_1},\overline{h_2},\overline{F}\right)=1 \mbox{ in } \F_q[x].\]
\end{thm}

\section{The Proof of Theorem \ref{Thm:Main1}}\label{Section:Main1Proof}

\begin{proof}
For each case, we let $K=\Q(\theta)$, where $f_t(\theta)=0$, and we let $\Z_K$ denote the ring of integers of $K$.
\subsection*{{\bf 4T2}}
Since
\[f_t(x-1)=x^4-4x^3+(4t+6)x^2-(8t+4)x+4t+2\]
is 2-Eisenstein, we conclude that $f_t(x)=x^4+4tx^2+1$ is irreducible over $\Q$. The cubic resolvent of $f_t(x)$ is
\[r_t(x)=x^3-4tx^2-4x+16t=(x-2)(x+2)(x-4t).\]
By Theorem \ref{Thm:KW1}, it follows that $\Gal(f_t)\simeq$ 4T2, which proves \eqref{X=2:I1}.

For item \eqref{X=2:I2}, a straightforward calculation reveals that
\begin{equation}\label{Eq:Delta(f)}
\Delta(f_t)=2^8(4t^2-1)^2.
\end{equation}
To establish the monogenicity of $f_t(x)$, we use Theorem \ref{Thm:Dedekind} with $T(x):=f_t(x)$ and we examine the prime divisors $q$ of $\Delta(f_t)$.

 First, let $q=2$. Then $\overline{T}(x)=(x+1)^4$ and we can let $h_1(x)=x+1$ and $h_2(x)=(x+1)^3$.
Then,
\[F(x)=\dfrac{(x+1)^4-f_t(x)}{2}=2x^3-(2t-3)x^2+2x \equiv x^2 \pmod{2}.\]
Thus, $\gcd(\overline{h_1},\overline{h_2},\overline{F})=1$ and $\left[\Z_K:\Z[\theta]\right]$ is not divisible by $q=2$.

Next, we give details only for the case when a prime $q$ divides $2t-1$ since the case when $q$ divides $2t+1$ is similar. Since $t\equiv 1/2 \pmod{q}$, it follows that
\[\overline{T}(x)=\left\{\begin{array}{cl}
  (x^2+1)^2 &\mbox{ if $q\equiv 3 \pmod{4}$}\\[.5em]
  (x-y)^2(x+y)^2 & \mbox{ if $q\equiv 1 \pmod{4}$,}
 \end{array}\right.\]
 where $y\in \Z$ such that $y^2\equiv -1\pmod{q}$. If $q\equiv 3 \pmod{4}$, we can let $h_1(x)=h_2(x)=x^2+1$. Then
   \[\overline{F}(x)=\overline{-2\left(\dfrac{2t-1}{q}\right)}x^2\] so that $\gcd(\overline{h_1},\overline{h_2},\overline{F})=1$ if and only if  $\left[\Z_K:\Z[\theta]\right]$ is not divisible by $q$ if and only if $q^2\nmid (2t-1)$. If $q\equiv 1 \pmod{4}$, then $y^2=zq-2t$ for some integer $z$ since $-2t\equiv -1 \pmod{q}$. Then,
    \[\overline{F}(x)=\overline{-2z}(x^2+1)+\overline{\left(\dfrac{4t^2-1}{q}\right)},\] from which we see that
    \[\overline{F}(\pm y)=\overline{\left(\dfrac{4t^2-1}{q}\right)}\ne 0 \mbox{ if and only if } q^2\nmid (2t-1).\] Thus, we conclude that $f_t(x)$ is monogenic if and only if $4t^2-1$ is squarefree, completing the proof of item \eqref{X=2:I2}.

For \eqref{X=2:I3}, we note first that $\FF_2$ is an infinite set since there exist infinitely many integers $t$ such that $4t^2-1$ is squarefree \cite{BB}. To see that the fields generated by the elements $f_t(x)$ of $\FF_2$ are all distinct, we proceed by way of contradiction. We assume for integers $t_1\ne t_2$, that $K_1=K_2$, where $K_1=\Q(\alpha_1)$ and $K_2=\Q(\alpha_2)$ with $f_{t_1}(\alpha_1)=0=f_{t_2}(\alpha_2)$. Since both $f_{t_1}(x)$ and $f_{t_2}(x)$ are monogenic by item \eqref{X=2:I2}, it follows that $\Delta(f_{t_1})=\Delta(f_{t_2})$. Consequently, from \eqref{Eq:Delta(f)}, we have that
\[8(t_1-t_2)(t_1+t_2)(2t_1^2+2t_2^2-1)=0.\] Hence, $t_1=-t_2$. Without loss of generality, we may assume that $t_1>0$ so that $t_2<0$. Since the zeros of $f_t(x)$ are
\[ \sqrt{-2t+\sqrt{4t^2-1}}, \sqrt{-2t-\sqrt{4t^2-1}}, -\sqrt{-2t+\sqrt{4t^2-1}}, -\sqrt{-2t-\sqrt{4t^2-1}},\]
it follows that all zeros of $f_t(x)$ are real if $t<0$, while all zeros of $f_t(x)$ are nonreal if $t\ge 0$, which contradicts the assumption that $K_1=K_2$.

\subsection*{{\bf 4T3}}
Note that
\[f_0(x-1)=x^4-4x^3+10x^2-8x+2\]
is 2-Eisenstein. Hence, $f_0(x)$ is irreducible over $\Q$. Suppose next that $t\ne 0$. Then $\abs{t}>10/12$, from which we deduce that
\[\abs{24t}=\abs{12t}+12\abs{t}>\abs{12t}+10\ge \abs{12t+4}+6.\] Thus, when $t\ne 0$, we have that $f_t(x)=x^4+24tx^3+(12t+4)x^2+4x+1$ is irreducible over $\Q$ by Perron's Irreducibility Criterion \cite{Perron}.

To determine $\Gal(f_t)$, we use Theorem \ref{Thm:KW1}. The cubic resolvent of $f_t(x)$ is
\begin{align*}
 r_t(x)&=x^3-(12t+4)x^2+(96t-4)x-(576t^2-48t)\\
 &=(x-12t)(x^2-4x+(48t-4)).
\end{align*}
Since
\[\Delta(x^2-4x+(48t-4))=-192t+32=16(-12t+2),\]
and $-12t+2 \equiv 2 \pmod{4}$, we conclude that $r_t(x)$ has exactly the one rational zero $x=12t$. Consequently, the splitting field of $r_t(x)$ is $L=\Q\left(\sqrt{-2(6t-1)}\right)$. Then, in Theorem \ref{Thm:KW1}, we have that $g(x)=g_1(x)g_2(x)$, where
\[g_1(x)=x^2-12tx+1 \mbox{ and } g_2(x)=x^2+24tx+4.\] Since $\Delta(g_1)=4(6t-1)(6t+1)$ and $\Delta(g_2)=16(6t-1)(6t+1)$, it follows that both $g_1(x)$ and $g_2(x)$ are irreducible over $L$. Hence, $\Gal(f_t)\simeq$ 4T3 by Theorem \ref{Thm:KW1}, which proves item \eqref{X=3:I1}.

To establish item \eqref{X=3:I2}, we use Theorem \ref{Thm:Dedekind} with $T(x):=f_t(x)$. A straightforward calculation yields
\[\Delta(f_t)=-2^9(6t-1)^3(6t+1)^2.\] First let $q=2$. Then $\overline{T}(x)=(x+1)^4$, and we can select
\begin{equation}\label{Eq:h}
h_1(x)=x+1 \mbox{ and } h_2=(x+1)^3.
\end{equation}
Thus,
\[F(x)=\dfrac{(x+1)^4-f_t(x)}{2}=(2-12t)x^3+(1-6t)x^2 \equiv x^2 \pmod{2}.\]
Therefore,
 $\gcd(\overline{h_1},\overline{F})=1$ so that $\left[\Z_K:\Z[\theta]\right]$ is not divisible by $q=2$.

 When $q$ is a prime divisor of $6t-1$, we also get that
 \[\overline{T}(x)=x^4+4x^3+6x^2+4x+1=(x+1)^4,\] with $h_i(x)$ as in \eqref{Eq:h}. Then,
 \[F(x)=\dfrac{(x+1)^4-f_t(x)}{2}=\dfrac{-2(6t-1)}{q}(2x+1)x^2,\] and it follows that
\[\gcd(\overline{h_1},\overline{F})=1 \mbox{ if and only if } q^2\nmid (6t-1).\] Thus,
$q\nmid \left[\Z_K:\Z[\theta]\right]$ if and only if $q^2\nmid (6t-1)$.

Suppose next that $q$ is a prime divisor of $6t+1$. Then $t\equiv -1/6 \pmod{q}$ and
\[T(x)\equiv (x^2-2x-1)^2 \pmod{q}.\] If $x^2-2x-1$ is irreducible over $\F_q$, then we can let $h_1(x)=h_2(x)=x^2-2x-1$. In this case, we get
\[ F(x)=\dfrac{(x^2-2x-1)^2-f_t(x)}{2}=-2\left(\dfrac{6t+1}{q}\right)x^2(2x+1),\]
so that
\[\gcd(\overline{h_1},\overline{F})=1 \mbox{ if and only if } q^2\nmid (6t+1).\]  If $x^2-2x-1$ is reducible over $\F_q$, then
\[x^2-2x-1=(x-(1+y))(x-(1-y)),\]
where $y\in \Z$ with $y^2\equiv 2\pmod{q}$. Then, we can select
\[h_1(x)=h_2(x)=(x-(1+y))(x-(1-y)),\] so that
\[F(x)=-4\left(\dfrac{6t+1}{q}\right)x^3-2\left(\dfrac{6t+1}{q}+\dfrac{y^2-2}{q}\right)x^2+4\left(\dfrac{y^2-2}{q}\right)x+y^2\left(\dfrac{y^2-2}{q}\right).\]
Then, computer calculations reveal that
\[\overline{F}(1\pm y)=2\left(\dfrac{6t+1}{q}\right)(\mp 12y-17).\] If $\mp 12y-17\equiv 0 \pmod{q}$, then, since $y^2\equiv 2\pmod{q}$, we get that
\[288\equiv 144y^2\equiv (\mp 12y)^2\equiv 17^2\equiv 289 \pmod{q},\] which yields the contradiction that $0\equiv 1\pmod{q}$. Hence,
\[\gcd(\overline{h_1},\overline{F})=1 \mbox{ if and only if } q^2\nmid (6t+1),\] completing the proof that
$f_t(x)$ is monogenic if and only if $36t^2-1$ is squarefree.

For \eqref{X=3:I3}, we proceed as in the case X=2. The set $\FF_3$ is infinite since there exist infinitely many values of $t$ such that $36t^2-1$ is squarefree \cite{BB}. We assume for integers $t_1\ne t_2$, that $K_1=K_2$, where $K_1=\Q(\alpha_1)$ and $K_2=\Q(\alpha_2)$ with $f_{t_1}(\alpha_1)=0=f_{t_2}(\alpha_2)$. Since both $f_{t_1}(x)$ and $f_{t_2}(x)$ are monogenic by item \eqref{X=3:I2}, it follows that $\Delta(f_{t_1})=\Delta(f_{t_2})$. Then, solving this discriminant equation using Maple shows that $6t_1$ must be a zero of the polynomial
\begin{multline*}
\GG(X):=X^4+(6t_2-1)X^3+(36t_2^2-6t_2-2)X^2\\+(216t_2^3-36t_2^2-12t_2+2)X
+1296t_2^4-216t_2^3+12t_2-72t_2^2+1,
\end{multline*} which is impossible since $\GG(0)\not \equiv 0 \pmod{6}$. This contradiction completes the proof of \eqref{X=3:I3} for the case of X=3.

\subsection*{{\bf 4T4}} Since $f_t(x)$ is 2-Eisenstein, $f_t(x)$ is irreducible over $\Q$. Straightforward calculations using Maple give us
\[\Delta(f_t)=2^6(4t-1)^2(108t^2-54t+7)^2,\] and the cubic resolvent  of $f_t(x)$ as
\[r_t(x)=x^3-2x^2-(144t^2-72t+8)x+128t^2-64t+8.\] If $r_t(x)$ is reducible over $\Q$, then
\[r_t(x)=(x+z)(x^2+Ax+B)=x^3+(A+z)x^2+(Az+B)x+Bz,\] for some integers $A,B,z$. Hence, by equating coefficients, we see that
\begin{align*}
  A+z&=-2\\
  Az+B&=-(144t^2-72t+8)\\
  Bz&=128t^2-64t+8.
\end{align*}
Solving this system in Maple tells us that
\[B=(8z^2+16z+8)/(9z+8).\] Let $d:=\gcd(8z^2+16z+8,9z+8)$. Then, easy gcd calculations show that $d\mid 8$. Since $B$ is an integer, we know that $d=9z+8$. Hence, the only possibilities for $z$ are $z=0$ and $z=-1$. However, checking the values of $-z$ in $r_t(x)$ gives $r(0)=8(4t-1)^2\ne 0$
 and $r(1)=-(4t-1)^2\ne 0$, since $t\in \Z$.
 Hence, $r_t(x)$ is irreducible over $\Q$, and it follows from Theorem \ref{Thm:KW1} that $\Gal(f_t)\simeq$ 4T4, which proves \eqref{X=4:I1}.

For item \eqref{X=4:I2}, we use Theorem \ref{Thm:Dedekind} with $T(x):=f_t(x)$ to check the prime divisors $q$ of $\Delta(f_t)$. 
Suppose first that $q=2$.
Then $\overline{T}(x)=x^4$, so that we can let $h_1(x)=x$ and $h_2(x)=x^3$, and we get that
\[F(x)=-x^3-x^2-2tx-18t^2+8t-1\equiv x^3+x^2+1 \pmod{2}.\] Hence, we see quite easily that $\gcd(\overline{h_1},\overline{F})=1$, from which we conclude that $\left[\Z_K:\Z[\theta]\right]$ is not divisible by $q=2$.

Next, let $q$ be a prime divisor of $4t-1$. Then, $t\equiv 1/4\pmod{q}$ and
\[T(x)\equiv (x^2+x+1/2)^2 \pmod{q}.\] Suppose first that $x^2+x+1/2$ is irreducible over $\F_q$. Note that in this situation, we must have $q\equiv 3 \pmod{4}$ since $-1$ is not a square. Then we can let
\[h_1(x)=h_2(x)=x^2+x+(q+1)/2,\] which yields
\begin{align*}
\overline{F}(x)&=x^2-\overline{\left(\dfrac{4t-1}{q}-1\right)}x-\overline{\left(\dfrac{(4t-1)(36t-7)-q^2-2q}{4q}\right)}\\
&=x^2-\overline{\left(\dfrac{4t-1}{q}-1\right)}x-\overline{\left(\dfrac{(4t-1)(36t-7)}{4q}\right)}+\overline{\dfrac{1}{2}}.
\end{align*} 
Note that $36t-7 \equiv 2\pmod{q}$ so that $q\nmid (36t-7)$. Hence, it follows that $\gcd(\overline{h_1},\overline{F})\ne 1$ if and only if $\overline{F}(x)=\overline{h_1}(x)$, which is true if and only if $q^2\mid (4t-1)$.

Suppose next that $x^2+x+1/2$ is reducible over $\F_q$. Then $q\equiv 1 \pmod{4}$ and
\[x^2+x+1/2\equiv\left(x-\left(\dfrac{-1+y}{2}\right)\right)\left(x-\left(\dfrac{-1-y}{2}\right)\right)\pmod{q},\]
where $y^2\equiv -1\pmod{q}$. Choosing $y\equiv 1 \pmod{2}$, we can let
\[h_1(x)=h_2(x)=\left(x-\left(\dfrac{-1+y}{2}\right)\right)\left(x-\left(\dfrac{-1-y}{2}\right)\right).\]
Then, computer calculations produce
\[\overline{F}(x)=\overline{\left(\dfrac{-y^2-1}{2q}\right)}x^2+\overline{\left(\dfrac{1-y^2-8t}{2q}\right)}x
+\overline{\left(\dfrac{y^4-2y^2-576t^2+256t-31}{16q}\right)}.\]
so that
\[\overline{F}((-1\pm y)/2)=\overline{\dfrac{-(4t-1)(9(4t-1)\pm 2y)}{q}+\dfrac{(y^2+1)^2}{4q}}=\overline{\mp 2y\left(\dfrac{4t-1}{q}\right)}.\]
Hence, we see that $\overline{F}((-1\pm y)/2)=0$ if and only if $q^2\mid (4t-1)$, which completes the proof when $q\mid (4t-1)$.

Now, suppose that $q$ is a prime divisor of $108t^2-54t+7$. Note that $q\not \in \{2,3\}$. Then, it follows that
\[\overline{T}(x)=(x-(18t-5))(x-(1-6t))^3,\] and therefore, we can let
\begin{equation}\label{Eq:h1h2}
h_1(x)=(x-(18t-5))(x-(1-6t)) \mbox{ and } h_2(x)=(x-(1-6t))^2.
\end{equation}
Observe from \eqref{Eq:h1h2}, that we only need to check $\overline{F}(1-6t)$ to determine $\gcd(\overline{h_1},\overline{h_2},\overline{F})$. More precisely,
\begin{equation}\label{Eq:Fbarcondition}
 \left[\Z_K:\Z[\theta]\right]\equiv 0 \Mod{q} \quad \mbox{if and only if} \quad \overline{F}(1-6t)=0.
\end{equation}
Straightforward calculations yield
\[F(x)=\left(\dfrac{108t^2-54t+7}{q}\right)\left(-2x^2-2(8t-1)x-(36t^2-10t+1)\right),\] so that
\[\overline{F}(1-6t)=-\overline{\left(\dfrac{108t^2-54t+7}{q}\right)}(12t^2-6t+1)=\overline{-\dfrac{2}{9}\left(\dfrac{108t^2-54t+7}{q}\right)},\] since \[9(12t^2-6t+1)-2=108^2-54t+7\equiv 0 \pmod{q}.\]
Hence,
\[\overline{F}(1-6t)=0 \quad \mbox{if and only if} \quad  108t^2-54t+7\equiv 0\Mod{q^2}.\] Consequently, from \eqref{Eq:Fbarcondition}, we have that
\[\left[\Z_K:\Z[\theta]\right]\equiv 0 \Mod{q} \quad \mbox{if and only if} \quad 108t^2-54t+7\equiv 0\Mod{q^2}.\] Since $\gcd(4t-1,108t^2-54t+7)=1$, the proof of item \eqref{X=4:I2} when X=4 is complete.

For item \eqref{X=3:I3}, we proceed as in the previous cases. The set $\FF_4$ is infinite since there exist infinitely many values of $t$ such that $(4t-1)(108t^2-54t+7)$ is squarefree \cite{BB}. We assume for integers $t_1\ne t_2$, that $K_1=K_2$, where $K_1=\Q(\alpha_1)$ and $K_2=\Q(\alpha_2)$ with $f_{t_1}(\alpha_1)=0=f_{t_2}(\alpha_2)$. Since both $f_{t_1}(x)$ and $f_{t_2}(x)$ are monogenic by item \eqref{X=3:I2}, it follows that $\Delta(f_{t_1})=\Delta(f_{t_2})$. Using Maple to solve this equation, we get six possible solutions. One solution has $t_1=t_2$, which we are not considering. A second solution has $t_1=1/2-t_2$, which is impossible in integers $t_1$ and $t_2$. The remaining four solutions contain the expression $\sqrt{-255+1944t_2-3888t_2^2}$. Thus, in order for one of these solutions to be viable, we must have that $-255+1944t_2-3888t_2^2$ is a perfect square. However, it is easy to see that $-3888x^2+1944x-255<0$ for all $x\in \R$, and
this fact completes the proof of \eqref{X=4:I3} for the case of X=4.

\subsection*{{\bf 4T5}}
For item \eqref{X=5:I1}, $f_t(x)$ is irreducible over $\Q$ for all $t\in \Z$ since $f_t(x)$ is 2-Eisenstein.
Using Maple, we calculate
\[\Delta(f_t)=16(4t+1)(4t-7)(64t+13).\] Since $\gcd(4t+1,4t-7)=1$ with $4t+1$ and $4t-7$ squarefree, it follows that $(4t+1)(4t-7)$ divides $64t+13$ if $\Delta(f_t)$ is a square. Observe then that $64t+13=16(4t+1)-3$ implies that $4t+1$ divides 3, from which we conclude that $t\in \{-1,0\}$. However, $64t+13=16(4t-7)+5^3$ implies that $4t-7$ divides $5$, which in turn implies that $t\in \{2,3\}$. This impossibility shows that $\Delta(f_t)$ is not a square. The cubic resolvent of $f_t(x)$ is
\[r_t(x)=x^3+2x^2-4(4t+1)x-12(4t+1).\]  If $r_t(x)$ is reducible over $\Q$, then
\[r_t(x)=(x+z)(x^2+Ax+B)=x^3+(A+z)x^2+(Az+B)x+Bz,\] for some integers $A,B,z$. Equating coefficients, we get that
\[A+z=2, \quad Az+B=-16t-4 \quad \mbox{and} \quad Bz=-48t-12.\] Thus,
\[B:=-\dfrac{3z(z-2)}{z-3}=-3(z+1)-\dfrac{9}{z-3}\in \Z,\]
which implies that $z-3$ divides 9, and hence, $-z\in \{-12,-6,-4,-2,0,6\}$. However, solving $r_t(-z)=0$ for $t$ for each of these values of $-z$ yields
\[t\in \{39/4,\ 11/4,\ 7/4,\ -1/4\},\] which contradicts the fact that $t\in \Z$. Therefore, $r_t(x)$ is irreducible for all $t\in \Z$, and we deduce from Theorem \ref{Thm:KW1} that $\Gal(f_t)\simeq$ 4T5.

To establish item \eqref{X=5:I2}, we use Theorem \ref{Thm:Dedekind} with $T(x):=f_t(x)$ to check the prime divisors $q$ of $\Delta(f_t)$. Since the details are similar for any prime divisor $q$ of $\Delta(f_t)$, we only give details in the case when $q$ divides $4t-7$. Then, $t\equiv 7/4 \pmod{q}$ and
\[T(x)\equiv (x^2-4x+5)(x+1)^2 \pmod{q}.\]  Since $x=-1$ is a zero of $x^2-4x+5$ if and only if $q=5$, we have three cases to consider:
\begin{equation}\label{Eq:4t-7}
\overline{T}(x)=\left\{\begin{array}{cl}
 x(x+1)^3 & \mbox{ if $q=5$}\\[.5em]
 (x^2-4x+5)(x+1)^2 & \mbox{ if $x^2-4x+5$ is irreducible over $\F_q$}\\[.5em]
 (x-A)(x-B)(x+1)^2 & \mbox{ if $x^2-4x+5$ is reducible over $\F_q$,}
\end{array}\right.
\end{equation}
where $A$ and $B$ are integers with $A,B\not \equiv -1 \pmod{q}$. We see in all of these cases that $h_2(x)=x+1$ so that  $\gcd(\overline{h_1},\overline{h_2},\overline{F})$ can be determined by calculating $\overline{F}(-1)$. Thus, more precisely, we have that
\[\gcd(\overline{h_1},\overline{h_2},\overline{F})=1\Longleftrightarrow \overline{F}(-1)\ne 0 \Longleftrightarrow \left[\Z_K:\Z[\theta]\right] \mbox{ is not divisible by $q$.}\]
Straightforward calculations reveal in all cases of \eqref{Eq:4t-7} that
\[\overline{F}(-1)=\overline{-\left(\dfrac{4t-7}{q}\right)}.\] Consequently, we deduce that
\[\left[\Z_K:\Z[\theta]\right]\not \equiv 0 \pmod{q} \quad \mbox{if and only if}\quad 4t-7 \mbox{ is squarefree}.\]

For \eqref{X=5:I3}, we first note that there exist infinitely many integers $t$ such that
 \[(4t+1)(4t-7)(64t+13)\] is squarefree \cite{BB}. Hence, there are infinitely many integers $t$ such that $4t+1$, $4t-7$ and $65t+13$ are simultaneously squarefree. Thus, the set $\FF_5$ is infinite. To see that each such value of $t$ generates a distinct field, we proceed as in previous cases, and let Maple solve the equation $\Delta(f_{t_1})=\Delta(f_{t_2})$. In all solutions given by Maple, other than $t_1=t_2$, we see that $W(t_2):=-12288t_2^2+10624t_2+19049$ must be a perfect square. It is easy to verify that $W(t_2)<0$ for all integers $t_2\not \in \{0,1\}$. Since neither $W(0)=19049$ nor $W(1)=17385$ is a square, the proofs of item \eqref{X=5:I3} and the theorem are complete.
\end{proof}

\section{Comparing Families in Theorem \ref{Thm:Main1} to Known Families}\label{Section:Known}
\subsection*{X=2} The literature on quartic polynomials with Galois group 4T2 is fairly extensive, with some authors addressing the issue of monogenicity  \cite{Chang,GT, JonesNYJM, Nyul, Williams}. However, we found only the following explicit two-parameter family of monogenic quartic 4T2-polynomials \cite{JonesNYJM} given by
\[\FF(x)=x^4+(36rp-1)x^2+1,\] where $r\ge 3$ and $p$ are primes, such that $r$ is a primitive root modulo 9 and
\[(12rp-1)(12rp+1)(36rp-1)(36rp+1)\]
is squarefree. Suppose that $\FF(\alpha)=0$ and $f_t(\theta)=0$. To see that $K\ne L$, where $K:=\Q(\alpha)$ and $L:=\Q(\theta)$, suppose to the contrary that $K=L$. Since both $\FF(x)$ and $f_t(x)$ are monogenic, it follows that
\[144(36rp+1)^2(12rp-1)^2=\Delta(\FF)=\Delta(K)=\Delta(f_t)=256(2t-1)^2(2t+1)^2,\]
which is clearly seen to be impossible by examining the power of 2 dividing each side. Hence, $K\ne L$.
\subsection*{X=3} Several authors have investigated integral bases for quartic 4T3-number fields \cite{JonesNYJM,Kable,HSW,Led}, but we found only the following (two-parameter) families of monogenic quartic 4T3-polynomials \cite{JonesNYJM} in the existing literature given by
\[\FF^+(x)=x^4+x^3+(100rp+1)x^2+x+1\quad \mbox{and} \quad \FF^-(x)=x^4-x^3+(100rp+1)x^2-x+1,\]
where $r\ge 3$ and $p$ are primes, such that $r$ is a primitive root modulo 25 and
\[(20rp+1)(100rp+1)(80rp-1)\]
is squarefree. Note that, for fixed values of $r$ and $p$, $\FF^+(x)$ and $\FF^-(x)$ generate the same field. Since
\[\Delta(\FF^{\pm})=5^3(20rp+1)(100rp+1)(80rp-1)^2 \quad \mbox{and} \quad \Delta(f_t)=-2^9(6t-1)^3(6t+1)^2,\] an argument similar to the one given for X=2 easily shows that the family $\FF^+(x)$ differs from the family given in Theorem \ref{Thm:Main1} for X=3.
\subsection*{X=4} 
The only monogenic 4T4-quartic family of polynomials we could find in the existing literature is the single-parameter family
\[\FF_m(x)=x^4+18x^2-4mx+m^2+81,\] where $m(m^2+108)$ is squarefree \cite{Spearman}. Note that $m$ is odd. To see that the family $\FF_m(x)$ has no overlap with the family presented in this article, we assume for some integers $m$ and $t$ for which $\FF_m(x)$ and $f_t(x)$ are respectively monogenic, that $K=L$, where $K:=\Q(\alpha)$ and $L:=\Q(\theta)$ with  $\FF_m(\alpha)=0$ and $f_t(\theta)=0$. Then
\[2^8m^2(m^2+108)^2=\Delta(\FF)=\Delta(K)=\Delta(f_t)=2^6(4t-1)^2(108t^2-54t+7)^2.\] Using Maple to solve this equation, we get that any solution must have that \begin{equation}\label{Eq:E1}
z^2=11664m^6+2519424m^4+136048896m^2+1,
\end{equation}
for some integer $z$. Multiplying both sides of \eqref{Eq:E1} by 4 shows that $(x,y)=(36m^2,2z)$ is an integral point on the elliptic curve $E_1$ given by \[y^2=x^3+7776x^2+15116544x+4.\] 
Using Sage to find all integral points $(x,y)$ on $E_1$, we get (with $y\ge 0$) that
\[(x,y)\in \{(0,2),(-3888,2),(14281868898720,53973124902433105922)\}.\] Since $x=36m^2$ with $m$ odd, it is easy to see that none of these points yields a valid solution to \eqref{Eq:E1}. That is, there is no overlap with these two families.
\subsection*{X=5}
We found two families of monogenic quartic 4T5-trinomials in \cite{Smith}, and one family in \cite{GSS}. The first family in \cite{Smith} is
\begin{equation}\label{Smith1}
\FF_{b}(x)=x^4+bx+b\in \Z[x] \quad \mbox{with} \quad \Delta(\FF_b)=(256-27b)b^3,
\end{equation} where $b\not \in \{3,5\}$, and both $b$ and $256-27b$ are squarefree. The second family in \cite{Smith} is
\begin{equation}\label{Smith2}
\GG_{d}(x)=x^4+x^3+d\in \Z[x] \quad \mbox{with} \quad \Delta(\GG_d)=(256d-27)d^2,
\end{equation} where $d\ne -2$, and both $d$ and $256d-27$ are squarefree. The family in \cite{GSS} is
\begin{equation}\label{GSS}
\HH_m(x)=x^4-6x^2-mx-3\in \Z[x] \quad \mbox{with} \quad \Delta(\HH_m)=-27(m-8)^2(m+8)^2.
\end{equation}
  Recall that $f_t(x):=x^4-2x^3-2x^2+6x+4t-2$, with \[\Delta(f_t)=16(4t+1)(4t-7)(64t+13).\]

  The question is whether $f_t(x)$ and $\FF_{b}(x)$ (or $\GG_d(x)$ or $\HH_m(x)$) generate the same quartic field for some integers $t$ and $b$ (or $d$ or $m$). If so, then $K=L$, where $K=\Q(\alpha)$, $L=\Q(\theta)$, $f_t(\theta)=0$ and $\FF_{b}(\alpha)=0$ (or $\GG_{d}(\alpha)=0$ or $\HH_m(\alpha)=0$). More importantly in this situation, since all of these polynomials are monogenic, we must have that their respective discriminants are equal. In each case, we assume that these polynomial discriminants are equal and proceed toward a contradiction.

  We begin with $\FF_b(x)$ in \eqref{Smith1}. Here we are assuming that
\begin{equation*}\label{Eq:DiscequationS1}
(256-27b)b^3=16(4t+1)(4t-7)(64t+13),
\end{equation*} for some integers $b$ and $t$. Since $4t+1$, $4t-7$ and $64t+13$ are all odd and squarefree, it follows that $b=2$. Then Maple tells us that $P(t)=0$, where
\[P(x):=128x^3-166x^2-95x-24.\] Since $\Delta(P)<0$, we know that $P(x)$ has exactly one real zero, and since $P(1)<0$ while $P(2)>0$, it follows that $P(t)$ has no integer zero, which contradicts the fact that $t\in \Z$.

For $\GG_d(x)$ in \eqref{Smith2}, a similar argument shows that $d=4$ and $P(t)=0$, where
\[P(x):=128x^3-166x^2-95x-136.\] As before, $P(x)$ has exactly one real non-integer zero between 1 and 2, which contradicts the fact that $t\in \Z$.

Finally, we address $\HH_m(x)$ in \eqref{GSS}. We assume that
\begin{equation}\label{Eq:Discequation1GSS}
-27(m-8)^2(m+8)^2=16(4t+1)(4t-7)(64t+13),
\end{equation} for some integers $m$ and $t$. An examination of each of the factors $A:=4t+1$, $B:=4t-7$ and $C:=64t+13$ modulo 3 
produces the following:
\[[A \Mod{3},B \Mod{3},C\Mod{3}]=\left\{\begin{array}{cl}
[1,2,1] & \mbox{if $t\equiv 0 \Mod{3}$}\\[.5em]
[2,0,2] & \mbox{if $t\equiv 1 \Mod{3}$}\\[.5em]
[0,1,0] & \mbox{if $t\equiv 2 \Mod{3}$.}
\end{array}\right.\] Thus, since $A$, $B$ and $C$ are squarefree, it follows that $ABC\not \equiv 0 \pmod{27}$, and so there are no solutions to equation \eqref{Eq:Discequation1GSS}.

\section{The Case 4T1}\label{Section:4T1}
The case of monogenic quartic polynomials $f(x)$ such that $\Gal(f)\simeq$ 4T1 seems to be quite different from the other possible quartic Galois groups. Gras \cite{GrasCyclic} showed that there are only two imaginary monogenic cyclic quartic fields
\[\Q(\zeta_5) \quad \mbox{and}\quad  \Q(\zeta_{16}-\zeta_{16}^{-1}),\] where $\zeta_n$ is a primitive $n$th root of unity. Respective corresponding monogenic polynomials are $\Phi_5(x)=x^4+x^3+x^2+x+1$ and $x^4+4x^2+2$.

Olajos \cite{Olajos} proved for the simplest quartics
\[\FF_k(x)=x^4-kx^3-6x^2+kx+1, \quad k\not \in \{\pm 3,0\},\] that there are only two values of $k$, namely $k\in \{2,4\}$, for which there exists a power integral basis for the field $\Q(\alpha)$, where $\FF_t(\alpha)=0$. Since these fields are real, they represent fields distinct from the two imaginary monogenic quartic fields found by Gras. Note that $\FF_2(x)$ and $\FF_4(x)$ are not monogenic. However, from \cite{Olajos}, we can easily construct two monogenic quartic polynomials $f_2(x)$ and $f_4(x)$ corresponding, respectively, to the two real monogenic cyclic quartic fields of Olajos:
\[f_2(x)=x^4-10x^3+25x^2-20x+5 \quad \mbox{and}\quad f_4(x)=x^4-8x^3+16x^2-8x-2.\] 
Easy calculations in Maple show that
\begin{equation*}\label{Eq:f2f4}
\Delta(f_2)=2^4\cdot 5^3 \quad \mbox{and} \quad \Delta(f_4)=2^{11}.
\end{equation*}  While we could not find a family of real monogenic 4T1-quartic polynomials, we did manage to find four additional distinct real monogenic quartics:
\begin{align*}
  g_1(x)&=x^4+9x^3+19x^2+9x+1 \quad \mbox{with} \quad \Delta(g_1)=3^2\cdot 13^3,\\
  g_2(x)&=x^4+5x^3+5x^2-5x-5 \quad \mbox{with} \quad \Delta(g_2)=3^2\cdot 5^3,\\
  g_3(x)&=x^4+11x^3+31x^2+11x+1 \quad \mbox{with} \quad \Delta(g_3)=5^3\cdot 11^2 \quad \mbox{and}\\
  g_4(x)&=x^4+7x^3+9x^2-7x+1 \quad \mbox{with} \quad \Delta(g_4)=5^3\cdot 7^2.
\end{align*}
Using Maple and Theorem \ref{Thm:Dedekind}, it is straightforward to verify that the polynomials $f_2(x)$, $f_4(x)$, $g_1(x)$, $g_2(x)$, $g_3(x)$ and $g_4(x)$ are all monogenic. By comparing discriminants and using Maple, it is also easy to see that these polynomials generate distinct real fields.

The authors of \cite{MNSU} have outlined two approaches to generate monogenic cyclic quartic fields, and they claim to prove that the number of such fields is infinite by providing an argument to show that the set of these fields from their second approach has positive density. However, according to our private communications with two analytic number theorists, Daniel White and Stanley Xiao, we believe their density argument is incorrect. In particular, we quote Stanley Xiao:
\begin{quote}
 Yes, as you suspect, their claim is wrong. One cannot simply compare the densities $\delta, \delta_1$ as they suggest, because in one set-up the input sequence is much thinner than the other. Indeed, for $d = d_1^2 + 4$ the input sequence has density $\sim N^{1/2}$, while for $d_1 = d_2^2 + 4$ the input sequence (with $d_2$ freely varying) has density $N$. They are not comparable.
\end{quote}
Consequently, we believe that the existence of a family of monogenic 4T1-quartic polynomials is still unresolved.





\end{document}